\documentclass[smallextended,natbib,runningheads]{article}
\usepackage{graphicx}

\usepackage{amsmath, amssymb}

\usepackage{subfigure}

\usepackage{float}
\usepackage{placeins}
\newcommand{\R}{\mathbb R}
\newcommand{\E}{\mathbb E}

\newtheorem{theorem}{\sc Theorem}[section]

\newtheorem{lemma}{\sc Lemma}[section]

\newtheorem{remark}{\sc Remark}[section]

 \def \cqfd{\hspace*{14.6cm}$\blacksquare$}
 
\begin{document}

\title{Nonparametric test for detecting change in distribution with
panel data}
\author{M. Boutahar,  B. Ghattas and D. Pommeret
\\
{\small {\sc Institute of Mathematics of Luminy.}}  \\
{\small {\sc  Luminy Faculty of Sciences.  163 Av. de Luminy 13288 Marseille  Cedex 9 - France}} \\{ \small boutahar@univmed.fr , ghattas@univmed.fr , pommeret@univmed.fr}}

\maketitle

\begin{abstract}
This paper considers the problem of comparing two processes with panel data.
A nonparametric test is proposed for detecting a
monotone change in the link between the two process distributions.
The  test statistic is of CUSUM type, based on the empirical distribution functions. The asymptotic distribution of the proposed statistic is derived and its finite sample property is examined by bootstrap procedures through Monte Carlo simulations.

\textbf{keywords}{nonparametric estimation \and panel data \and process}

\end{abstract}

\section{Introduction}

Many situations lead to the comparison of two random processes.
In a parametric case, the problem of change detection has been widely
studied in the time series literature. A common problem is to test a change in
the mean or in the variance of the time series by using a parametric model
(see for instance \cite{gom} or \cite{gal}, and references
therein). In the Gaussian case comparisons
of processes are considered  through their covariance
structures (see \cite{gup}, \cite{pan}).
These distribution assumptions can be relaxed when the study
concerns processes observed through  panel data. This situation is
frequently encountered in medical follow-up studies when two groups of
patients are observed and compared. Each subject in the study gives rise to
a random process $(X_t)$ denoting the measurement of the patient up to time $%
t$ (such data are referred to as panel data).
In this context, \cite{bal8,bal9,bal10} considered the problem of testing the equality of mean functions and proposed  new multi-sample tests for panel count data.

In this paper we consider the general problem of comparison of two processes
which may differ by a transformation of their distributions. Our purpose is to test
whether this  transformation changes over time. For this, two panels are considered: $%
(X_{i,t})_{1\leq i\leq N_{x};1\leq t\leq n}$ and $(Y_{i,t})_{1\leq i\leq
N_{y};1\leq t\leq n}$, not necessarily independent; that is, we can have
i.i.d. paired observations $(X_i,Y_i)_{i=1,\cdots, N}$ with dependence
between $X_i$ and $Y_i$. It is  assumed that for each $t$, the $X_{i,t},1\leq
i\leq N_{x}$ (resp. $Y_{i,t},1\leq i\leq N_{y})$ are i.i.d. random variables
with common distribution function $F_{t}$ (resp. $G_{t}$) and with support ${\cal X}$ (resp. ${\cal Y}$). Also we assume
that for all $1\leq t\leq n$ there exists monotone transformations $h_{t}$ such that the following
equality in distribution holds: $X_{t}=^{\!\!\!\!d} h_{t}(Y_{t})$.
Without loss of generality we consider that the functions $h_t(.)$ are  increasing. Note
that if $F_{t}$ is invertible then there exists a trivial transformation $%
h_{t}$ given by $h_{t}=F_{t}^{-1}\circ G_{t}$.
We are interested in testing
whenever this transformation is time independent; that is, for all $t$, the equality $h_t=h$ occurs. A simple illustration  is the case where $X_t$ and $Y_t$ are  Gaussian processes with mean $m_X$ and $ m_Y$
and variance $t\sigma_X^2$ and $t \sigma_Y^2$, respectively. In that case the function $h$ is linear.

More generally, observing both processes $X_{t}$ and $Y_{t}$ with panel data we want to test
\begin{eqnarray*}
H_{0}:\forall t,h_{t}=h & \mathrm{against } & H_{1}:\exists t_{1}\neq
t_{2},h_{t_{1}}\neq h_{t_{2}}.
\end{eqnarray*}
It is clear that $H_0$ coincides with the equality in distribution:  $X_{t}=^{\!\!\!\!d}(h(Y_{t}))$, for all $t$.
Following \cite{gom} (see also \cite{gal}), we construct a non parametric  test statistic  based on the empirical
estimator of $h_t$, denoted by $\widehat{h}_t$.  We show that $\widehat{h}_t$ is proportional to a Brownian bridge under $H_0$.

When $H_0$ is not rejected, it is of
interest to estimate $h$ and to interpret its estimator $\widehat{h}$. Then
this test  can be viewed  as a first step permitting to legitimate
estimation and interpretation of a constant transformation $h$ between the distributions
of two samples, possibly paired.


The paper is organized as follows: In Section 2 we construct the test
statistic. In Section 3 we perform a simulation study
using a bootstrap procedure to evaluate the finite sample property of the test.
The power is evaluated against alternatives where there are smooth scale or
position time changes in the process distribution. Section 4 contains brief concluding remarks.

\section{The test statistic}

A natural nonparametric estimator of $h_{t}$ is given by
\begin{eqnarray*}
\widehat{h}_{t}(\cdot )&=&X_{(N_{x}\widehat{G}_{t}(\cdot )),t},
\end{eqnarray*}%
where $X_{(i),t}$ denotes the $i$th order statistic and $\widehat{G}_{t}$
is the empirical distribution function of $(Y_{i,t})_{1\leq i\leq N_{y}}$, that is
\begin{eqnarray*}
\widehat{G}_{t}(x)&=&\frac{1}{N_{y}}\sum_{i=1}^{N_{y}}\mathbf{1}%
_{\{Y_{i,t}\leq x\}}.
\end{eqnarray*}
A nonparametric test is considered to test the variation of $h_t$. For $\tau \in (0,1)$, $x\in \mathcal{Y}$,
write
\begin{eqnarray}
B_{n}(\tau ,x)&=&\frac{1}{\sqrt{n}\widehat{\sigma }_{n}}\left(
\sum_{t=1}^{[n\tau ]}\widehat{h}_{t}(x)-\frac{[n\tau ]}{n}\sum_{t=1}^{n}%
\widehat{h}_{t}(x)\right) ,  \label{bn}
\end{eqnarray}%
where
\begin{eqnarray*}
\widehat{\sigma }_{n}^{2}&=&\frac{1}{n}\sum_{t=1}^{n}(\widehat{h}_{t}(x)-\bar{h%
}(x))^{2},\bar{h}(x)=\frac{1}{n}\sum_{t=1}^{n}\widehat{h}_{t}(x).
\end{eqnarray*}

For a given square integrable function $w$ we define the following test
statistic

\begin{eqnarray*}
S_{n}(w)&=&\int_{
\R
}w(x)\sup_{1\leq \tau \leq 1}\left\vert B_{n}(\tau ,x)\right\vert dx.
\end{eqnarray*}%
To establish the limiting distribution of the statistic $S_{n}(w)$ under the
null, we need the following assumptions:
\begin{itemize}
\item
Assumption 1. There exists $a< \infty$ such that $N_{x}/(N_{x}+N_{y})\rightarrow a.$
\item
Assumption 2. There exist $\gamma _{1}>0$ and $\gamma _{2}>0$ such that $%
f_{t}(x)\geq \gamma _{1}$ and $g_{t}(y)\geq \gamma _{2}$ for all $(x,y) \in {\cal X} \times {\cal Y}$, where $f_{t}$ and $g_{t}$ are
the density \ functions of $X_{t}$ and $Y_{t}.$
\item
Assumption 3. 
For all $x\in {\cal X}$, there exist $ 0< \overline{\sigma }
_{2}^{2}(x) < \infty $  such that
\begin{eqnarray*}
\ \frac{1}{n}\sum_{t=1}^{n}\sigma _{1,t}^{2}(x)\rightarrow \overline{\sigma }%
_{2}^{2}(x),&&\text{ as }n\rightarrow \infty,  
\end{eqnarray*}
where
\begin{eqnarray}
\sigma _{1,t}^{2}(x)&=&\sigma _{t}^{2}(x)\frac{N_{x}+N_{y}}{N_{x}N_{y}},\text{
and }\sigma _{t}^{2}(x)=\frac{G_{t}(x)(1-G_{t}(x))}{f_{t}^{2}(h_{t}(x))}.
\label{sig}
\end{eqnarray}
\item
Assumption 4.
\begin{eqnarray*}
\frac{n(N_{x}+N_{y})}{N_{x}N_{y}}&\rightarrow& 0.
\end{eqnarray*}
\end{itemize}
\begin{remark}
Assumptions 1 and 2 are standard. 
Assumption 3 states that the second moments  converge on average. 
 If Assumption 1 is satisfied, Assumption 4 is equivalent to $n=o(N_x)$ or $n=o(N_y)$.
\end{remark}
\begin{theorem}
\label{theo}Let assumptions 1-4 hold. Then under the null $H_{0}$ we
have the following convergence in distribution
\begin{eqnarray}
S_{n}(w)\rightarrow ^{\!\!\!\!\!d}S(w)&=&B_{\infty }\int_{
\mathbb{R}
}w(x)dx,\text{as }n\rightarrow \infty ,N_{x}\rightarrow \infty \text{ and }%
N_{y}\rightarrow \infty ,
\label{snf}
\end{eqnarray}%
where $B_{\infty }=\sup_{0\leq \tau \leq 1}|B(\tau )|$, and $B$ is a Brownian
bridge.
\end{theorem}
\begin{remark}
The cumulative distribution function of $B_{\infty }$ is given by
(see \cite{bil})
\begin{eqnarray*}
F_{B_{\infty }}(z)&=&1+2\sum_{k=1}^{\infty }(-1)^{k}\exp \{-2k^{2}z^{2}\}.
\end{eqnarray*}
\end{remark}
Before proving Theorem 1, we state  three lemmas.
\begin{lemma}
\label{lem1} Under Assumption 1 we have
\begin{eqnarray}
\left( \frac{N_{x}N_{y}}{N_{x}+N_{y}}\right) ^{1/2}(\widehat{h}%
_{t}(x)-h_{t}(x))\rightarrow ^{\!\!\!\!\!d}N(0,\sigma _{t}^{2}(x)),&&\text{ as
}N_{x}\rightarrow \infty ,N_{y}\rightarrow \infty ,  \label{nor}
\end{eqnarray}
where $\sigma _{t}^{2}(x)$ is given by (\ref{sig}).
\end{lemma}
\paragraph{Proof}  $(\mathbf{1}_{\{Y_{i,t}\leq x\}})$ is an i.i.d sequence with mean $%
G_{t}(x)$ and variance $G_{t}(x)(1-G_{t}(x)),$  hence an immediate
application of the central limit theorem yields
\begin{eqnarray}
N_{y}^{1/2}\left( \widehat{G}_{t}(x)-G_{t}(x)\right) &\rightarrow
^{\!\!\!\!\!d}&N(0,G_{t}(x)(1-G_{t}(x))).  \label{co1}
\end{eqnarray}
By the delta-method the last convergence implies that
\begin{eqnarray}
N_{y}^{1/2}\left( F^{-1}\left( \widehat{G}_{t}(x)\right) -F^{-1}\left(
G_{t}(x)\right) \right) &\rightarrow ^{\!\!\!\!\!d}&N(0,\sigma _{t}^{2}(x)).
\label{co2}
\end{eqnarray}
For $\ p\in]0;1[$ fixed, denote by $\widehat{F}_{t}^{-1}(p)$ the
sample $p$-quantile; that is,
$\widehat{F}_{t}^{-1}(p)=X_{(r),t}$, where $r=\left[ N_{x}p\right] +1$.
By Theorem 3 of \cite{sen} we obtain
\begin{eqnarray}
N_{x}^{1/2}(\widehat{F}_{t}^{-1}(p)-F_{t}^{-1}(p))\rightarrow
^{\!\!\!\!\!d}N\left( 0,\frac{p(1-p)}{f_{t}^{2}(F_{t}^{-1}(p))}\right) ,&&
\forall p\in (0,1).  \label{co3}
\end{eqnarray}
Let $\phi _{X}(t)=\E(\exp (itX))$
denotes the characteristic  function
of the random variable $X$ and let
  $\phi _{X|Y}(t)=\E(\exp (itX)\mid Y)$
denotes the conditional characteristic function
of the random variable $X$ conditional on $Y$. We have
\begin{eqnarray*}
\widetilde{H}_{t} &=&\left( \frac{N_{x}N_{y}}{N_{x}+N_{y}}\right) ^{1/2}(%
\widehat{h}_{t}(x)-h_{t}(x)) \\
&=&\left( \frac{N_{x}N_{y}}{N_{x}+N_{y}}\right) ^{1/2}\left( \widehat{F}%
_{t}^{-1}(\widehat{G}_{t}(x))-F_{t}^{-1}(G_{t}(x))\right) \\
&=&\widetilde{H}_{1,t}+\widetilde{H}_{2,t},
\end{eqnarray*}
where
\begin{eqnarray*}
\widetilde{H}_{1,t} &=&\left( \frac{N_{y}}{N_{x}+N_{y}}\right)
^{1/2}N_{x}^{1/2}\left( \widehat{F}_{t}^{-1}(\widehat{G}_{t}(x))-F_{t}^{-1}(%
\widehat{G}_{t}(x))\right)
\end{eqnarray*}
\begin{eqnarray*}
\widetilde{H}_{2,t} &=&\left( \frac{N_{x}}{N_{x}+N_{y}}\right)
^{1/2}N_{y}^{1/2}\left( F_{t}^{-1}(\widehat{G}_{t}(x))-F_{t}^{-1}(G_{t}(x))%
\right) .
\end{eqnarray*}%
Then we get
\begin{eqnarray*}
\phi _{\widetilde{H}_{t}}(u) &=&\E(\exp (iu\widetilde{H}_{t})) \\
&=&\E\left( \E\left[ \exp (iu\widetilde{H}_{t})\mid Y_{t}\right] \right) \\
&=&\E\left( \exp (iu\widetilde{H}_{2,t}\ )\ \E\left[ \exp (iu\widetilde{H}%
_{1,t}\ )\mid Y_{t}\right] \right) .
\end{eqnarray*}%
Moreover
\begin{eqnarray}
\E\left[ \exp (iu\widetilde{H}_{1,t})\mid Y_{t}\right] &=&\phi _{\widetilde{H}%
_{1,t\mid Y_{t}}}(u)  \label{f1} \\
&=&\phi _{N_{x}^{1/2}\left( \widehat{F}_{t}^{-1}(\widehat{G}%
_{t}(x))-F_{t}^{-1}(\widehat{G}_{t}(x))\right) \mid Y_{t}}\left( \left(
N_{y}/(N_{x}+N_{y})\right) ^{1/2}u\right)  \notag
\end{eqnarray}%
From (\ref{co3}) it follows that, $\forall v \in \R$,
\begin{eqnarray}
\phi _{N_{x}^{1/2}\left( \widehat{F}_{t}^{-1}(\widehat{G}_{t}(x))-F_{t}^{-1}(%
\widehat{G}_{t}(x))\right) \mid Y_{t}}(v)&\longrightarrow &
 \exp \left( -\frac{1}{2}v^{2}\widehat{\sigma }_{t}^{2}(x)\right),
\label{f2}
\end{eqnarray}%
as $N_{x}\rightarrow \infty$,   where%
\begin{eqnarray*}
\widehat{\sigma }_{t}^{2}(x)&=&\frac{\widehat{G}_{t}(x)(1-\widehat{G}_{t}(x))}{%
f_{t}^{2}\left( F_{t}^{-1}(\widehat{G}_{t}(x))\right) }.
\end{eqnarray*}%
The convergence (\ref{co1}) yields  $\widehat{G}_{t}(x)$ $\
\underrightarrow{\text{ }P\text{ \ }}$ $\ G_{t}(x)$, as $N_{y}\rightarrow
\infty ,$ which implies, combined with (\ref{f1})-(\ref{f2}), Assumption 1
and $h_{t}(x)=F_{t}^{-1}(G_{t}(x))$, that
\begin{eqnarray}
\E\left[ \exp (iu\widetilde{H}_{1,t})\mid Y_{t}\right]
&\underrightarrow{d} & \exp \left( -\frac{1}{2}%
(1-a)u^{2}\sigma _{t}^{2}(x)\right),   \label{co4}
\end{eqnarray}
as $N_{x}\rightarrow \infty$
and $N_{y}\rightarrow \infty$.
Moreover we have
\begin{eqnarray*}
\exp (iu\widetilde{H}_{2,t}\ )&=&\exp \left[ iu\left( \frac{N_{x}}{N_{x}+N_{y}}%
\right) ^{1/2}N_{y}^{1/2}\left( F_{t}^{-1}(\widehat{G}%
_{t}(x))-F_{t}^{-1}(G_{t}(x))\right) \right].
\end{eqnarray*}%
Since the function $x \mapsto \exp (iu x )$ is continuous, then the convergence (%
\ref{co2}) and Assumption 1 yield
\begin{eqnarray}
\exp (iu\widetilde{H}_{2,t}\ )\rightarrow ^{\!\!\!\!\!d}\ \exp
(iua^{1/2}H_{2,t}\ ),&&\text{as }N_{x}\rightarrow \infty ,N_{y}\rightarrow
\infty ,  \label{co5}
\end{eqnarray}%
where $H_{2,t}$ is centered Gaussian distributed  with variance equal to $\sigma _{t}^{2}(x))$.
From (\ref{co4})  and (\ref{co5}) it follows that, as $N_{x}\rightarrow
\infty$ and  $N_{y}\rightarrow \infty$,
\begin{eqnarray}
\exp (iu\widetilde{H}_{2,t}\ )\E\left[ \exp (iu\widetilde{H}_{1,t})\mid Y_{t}%
\right]  && \nonumber \\
\rightarrow ^{\!\!\!\!\!d} \ \exp (iua^{1/2}H_{2,t}\ )&&\exp \left( -%
\frac{1}{2}(1-a)u^{2}\sigma _{t}^{2}(x)\right) .  \label{co6} \\
&&  \notag
\end{eqnarray}%
Since $\E\left[ \exp (iu\widetilde{H}_{1,t})\mid Y_{t}\right] $ and $\exp (iu%
\widetilde{H}_{2,t})$ are bounded almost surely, it follows from (\ref{co6})
that
\begin{eqnarray*}
\phi _{\widetilde{H}_{t}}(u) &=&\E\left( \exp (iu\widetilde{H}_{2,t}\ )\E\left[
\exp (iu\widetilde{H}_{1,t})\mid Y_{t}\right] \right) \\
&\rightarrow &\E\left( \exp \left( iua^{1/2}H_{2,t}\right) \exp \left( -\frac{%
1}{2}(1-a)u^{2}\sigma _{t}^{2}(x)\right) \right) \text{, as }%
N_{x}\rightarrow \infty ,N_{y}\rightarrow \infty \\
&=&\exp \left( -\frac{1}{2}au^{2}\sigma _{t}^{2}(x)\right) \exp \left( -%
\frac{1}{2}(1-a)u^{2}\sigma _{t}^{2}(x)\right) \\
&=&\exp \left( -\frac{1}{2}u^{2}\sigma _{t}^{2}(x)\right),
\end{eqnarray*}%
therefore the desired conclusion (\ref{nor}) holds.
\\
\cqfd

Lemma \ref{lem1} implies that
\begin{eqnarray}
\widehat{h}_{t}(x)&=&h_{t}(x)+\sigma _{1,t}(x)\varepsilon _{t}+r_{t},
\label{mo1}
\end{eqnarray}
where $\sigma _{1,t}^{2}(x)$ is given by (\ref{sig}), $(\varepsilon
_{t})$ is a standard Gaussian white noise and the remainder term $r_{t}$ is
such that
\begin{eqnarray}
r_{t}&=&O_{P}\left( \left\{ (N_{x}+N_{x})/N_{x}N_{y}\right\} ^{1/2}\right) .
\label{rt}
\end{eqnarray}
Let $D=D[0,1]$\ be the space of random functions that are right-continuous
and have left limits, endowed with the Skorohod topology. The weak
convergence of a sequence of random elements $X_{n}$ in $D$ to a random
element $X$ in $D$ will be denoted by $X_{n}\Longrightarrow X.$ Let
\begin{eqnarray}
W_{n}(\tau )&=&\frac{1}{\ \sqrt{n}}\displaystyle \sum\limits_{t=1}^{[n\tau ]}\sigma
_{1,t}(x)\varepsilon _{t},\text{ \ \ \ }\tau \in \lbrack 0,1].  \label{xy}
\end{eqnarray}

\begin{lemma}
\label{lem2}
Under Assumptions 1-3 we have
\begin{eqnarray}
W_{n}&\Longrightarrow&  \overline{\sigma }_{2}(x)\ W,  \label{wn}
\end{eqnarray}%
where $W$ stands for the standard Brownian motion.
\end{lemma}
\paragraph{Proof} Assumption 2 implies that
\begin{eqnarray*}
\sigma _{1,t}^{2}(x) &\leq &\frac{1}{\gamma _{1}^{2}}\frac{N_{x}+N_{y}}{%
N_{x}N_{y}} \\
&\leq &C,
\end{eqnarray*}%
for some positive constant $C$ and $N_{x}$ and $N_{y}$ large enough. \ Hence
$\sigma _{1,t}^{2}(x)$ is a bounded deterministic sequence, therefore the
weak convergence (\ref{wn}) follows from Theorem A.1 of \cite{bou}.
\\
\cqfd
\begin{lemma}
\label{lem3}
Under the null $H_{0}$, as $n\rightarrow \infty$, $N_{x}\rightarrow \infty$ and
$N_{y}\rightarrow \infty$ we have
\begin{eqnarray}
\widehat{\sigma }_{n}^{2}=\frac{1}{n}\sum_{t=1}^{n}(\widehat{h}_{t}(x)-\bar{h%
}(x))^{2} & \underrightarrow{\text{ }d\text{ \ }}& \overline{\sigma }%
_{2}^{2}(x).  \label{csig}
\end{eqnarray}
\end{lemma}
\paragraph{Proof}
Under the null $H_{0}$: $h_{t}(x)=h(x)$ the equality (\ref%
{mo1}) becomes
\begin{eqnarray*}
\widehat{h}_{t}(x)&=&h(x)+\sigma _{1,t}(x)\varepsilon _{t}+r_{t}.
\end{eqnarray*}
Let $y_{t}=h(x)+\sigma _{1,t}(x)\varepsilon _{t}$ , \ $\overline{y}%
=\sum_{t=1}^{n}y_{t}/n,$ then by using the same argument as in Theorem 1 of
\cite{bou} we obtain%
\begin{eqnarray}
\frac{1}{n}\sum_{t=1}^{n}(y_{t}-\overline{y})^{2}&\underrightarrow{d}&%
 \overline{\sigma }_{2}^{2}(x).  \label{cvy}
\end{eqnarray}
We have
\begin{eqnarray}
\frac{1}{n}\sum_{t=1}^{n}(\widehat{h}_{t}(x)-\bar{h}(x))^{2}&& \nonumber \\
=\frac{1}{n}%
\sum_{t=1}^{n}(y_{t}-\overline{y})^{2}&+&\frac{1}{n}\sum_{t=1}^{n}(r_{t}-%
\overline{r})^{2}+2\frac{1}{n}\sum_{t=1}^{n}(y_{t}-\overline{y})(r_{t}-%
\overline{r}),  \label{som}
\end{eqnarray}
where $\overline{r}=\sum_{t=1}^{n}r_{t}/n.$ From (\ref{rt}) it follows  that
\begin{eqnarray*}
\overline{r} &=&O_{P}(\left( (N_{x}+N_{y})/N_{x}N_{y}\right) ^{1/2}) \\
&=&o_{p}(1), \text{ as }\ N_{x}\rightarrow \infty ,N_{y}\rightarrow \infty ,
\end{eqnarray*}
which implies that
\begin{eqnarray}
\frac{1}{n}\sum_{t=1}^{n}(r_{t}-\overline{r})^{2}=o_{p}(1),&&\text{ as }\
N_{x}\rightarrow \infty ,N_{y}\rightarrow \infty .  \label{rt2}
\end{eqnarray}
By using the Cauchy Shwartz inequality, we have
\begin{eqnarray*}
\frac{1}{n}\sum_{t=1}^{n}(y_{t}-\overline{y})(r_{t}-\overline{r})&\leq& \left(
\frac{1}{n}\sum_{t=1}^{n}(y_{t}-\overline{y})^{2}\right) ^{1/2}\left( \frac{1%
}{n}\sum_{t=1}^{n}(r_{t}-\overline{r})^{2}\right) ^{1/2}.
\end{eqnarray*}
Hence by  using  (\ref{cvy}) and (\ref{rt2}) we get
\begin{eqnarray}
\frac{1}{n}\sum_{t=1}^{n}(y_{t}-\overline{y})(r_{t}-\overline{r})=o_{p}(1),%
&&\text{ as }\ N_{x}\rightarrow \infty ,N_{y}\rightarrow \infty .  \label{yr}
\end{eqnarray}
The desired conclusion (\ref{csig}) holds by combining (\ref{cvy})-(\ref{yr}%
).
\\
\cqfd

\paragraph{Proof of Theorem 1} Under the null, the process $B_{n}(\tau ,x)$ in
(\ref{bn})  can be rewritten as
\begin{eqnarray*}
B_{n}(\tau ,x) &=&\frac{1}{\sqrt{n}\widehat{\sigma }_{n}}\left(
\sum_{t=1}^{[n\tau ]}\sigma _{1,t}(x)\varepsilon _{t}-\frac{[n\tau ]}{n}%
\sum_{t=1}^{n}\sigma _{1,t}(x)\varepsilon _{t}\right) +R_{n}(\tau ,x) \\
&=&\frac{1}{\widehat{\sigma }_{n}}\left( W_{n}(\tau )-\frac{[n\tau ]}{n}%
W_{n}(1)\right) +R_{n}(\tau ,x),
\end{eqnarray*}
where the remainder term $R_{n}(\tau ,x)$ is given by
\begin{eqnarray*}
R_{n}(\tau ,x)&=&\frac{1}{\sqrt{n}\widehat{\sigma }_{n}}\left(
\sum_{t=1}^{[n\tau ]}r_{t}-\frac{[n\tau ]}{n}\sum_{t=1}^{n}r_{t}\right).
\end{eqnarray*}
Now observe that
\begin{eqnarray*}
\sum_{t=1}^{[n\tau ]}r_{t}&=&O_{P}\left( [n\tau ]\left(
(N_{X}+N_{Y})/N_{X}N_{Y}\right) ^{1/2}\right) ,
\end{eqnarray*}%
which together with (\ref{csig}) implies that%
\begin{eqnarray*}
R_{n}(\tau ,x) &=&O_{P}\left( \left\{ n(N_{x}+N_{y})/N_{x}N_{y}\right\}
^{1/2}\right) , \\
&=&o_{p}(1)\text{ under assumption 4.}
\end{eqnarray*}
Hence
\begin{eqnarray*}
R_{n}(\tau ,x)&=&\frac{1}{\widehat{\sigma }_{n}}\left( W_{n}(\tau )-\frac{%
[n\tau ]}{n}W_{n}(1)\right) +o_{p}(1),
\end{eqnarray*}
which combined  with (\ref{wn}) and (\ref{csig}) yields
\begin{eqnarray*}
B_{n}(.,x)&\Longrightarrow&  B,
\end{eqnarray*}%
where $B(\tau )=W(\tau )- \tau  W(1)$ is a Brownian bridge.
 Therefore
\begin{eqnarray}
\sup_{1\leq \tau \leq 1}\left\vert B_{n}(\tau ,x)\right\vert &\rightarrow
^{\!\!\!\!\!d}&\sup_{1\leq \tau \leq 1}\left\vert B(\tau )\right\vert .
\label{cbn}
\end{eqnarray}%
Let $F(
\mathbb{R},
\mathbb{R}
)$ be the space of square integrable functions endowed with the uniform norm $%
\left\Vert .\right\Vert _{\infty }.$ For a given square integrable function $%
w$, the functional $\mathcal{G}_{w}$: $\left( F(
\mathbb{R},
\mathbb{R}
),\left\Vert .\right\Vert _{\infty }\right) \rightarrow (%
\mathbb{R}
,\left\vert .\right\vert )$ defined by%
\begin{eqnarray*}
\mathcal{G}_{w}(g)&=&\int_{%
\mathbb{R}
}w(x)g(x)dx,
\end{eqnarray*}
 is continuous. To obtain the convergence (\ref{snf}) it is
sufficient to apply (\ref{cbn}) and the continuous mapping theorem.
\\
\cqfd
\section{Empirical study}

For simplicity we consider $N_{x}=N_{y}=N$.  Data are generated from three models: first, $Y_t$ is  normally distributed with mean $0$ and variance $1$, and  $X_t$ is generated independently by the transformation $X_t=h_t(Z_t)$, where $Z_t$ is another Gaussian process with mean $0$ and variance $1$. Second, $Y_t$ is an autoregressive process of order 1 (AR1) with correlation  coefficient equal to 0.5, and  $X_t$ is generated independently by the transformation $X_t=h_t(Z_t)$, where $Z_t$ is another AR1 process. For the last model random variables are paired: $Y_t$ are independent  Gaussian variables with mean $0$ and variance $1$, and  $X_t=h_t(Y_t)$, that is, the time transformation is on the random variables. It is clear that this implies the same transformation for the corresponding distributions. 

\paragraph{Alternatives.}
The following five alternatives are considered

~\\
{\bf First alternative: A1}\\ Change in the mean.
$h_{1,t}(x)  =   \displaystyle\frac{2t^2}{1+t^2} + x$.

~\\
{\bf Second alternative: A2}\\ Change in the variance.
$h_{2,t}(x)  =   \displaystyle\frac{2t^2}{1+t^2} \; x$.

~\\
{\bf Third alternative: A3 } \\ Jump.
$h_{3,t}(x)  =   x + 0.05 t \mathbb I_{t<n/2} + 0.005(n-t)\mathbb I_{t\geq n/2}$, \\where $\mathbb I_{t\geq n/2}=1$ if ${t\geq n/2}$ and 0 otherwise.

~\\
{\bf Fourth alternative: A4} \\
Smooth change in the mean.
$ h_{4,t}(x) = x+(1+\exp(-0.01(t-1)))^{-1}$

~\\
{\bf Fifth alternative: A5} \\
Smooth change in the mean.
$ h_{5,t}(x) = x+(1+\exp(-0.05(t-1)))^{-1}$



\begin{figure}[H]
	\begin{center}
\subfigure[Time length = 20]
{
		\includegraphics[scale=0.7]{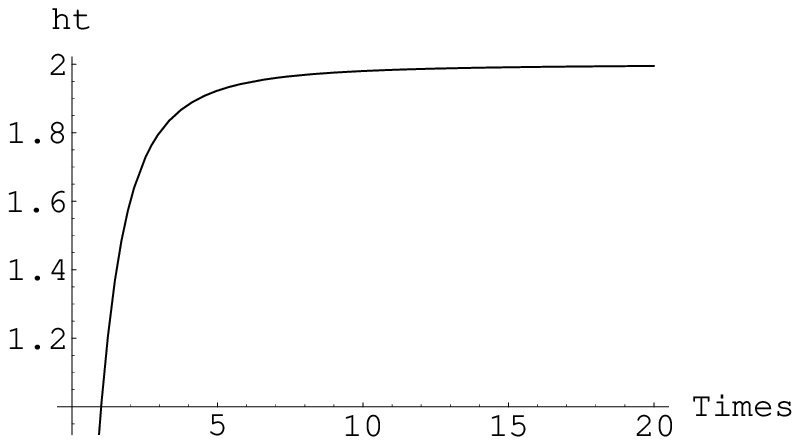}
}
%
\subfigure[Time length = 200]
{		\includegraphics[scale=0.68]{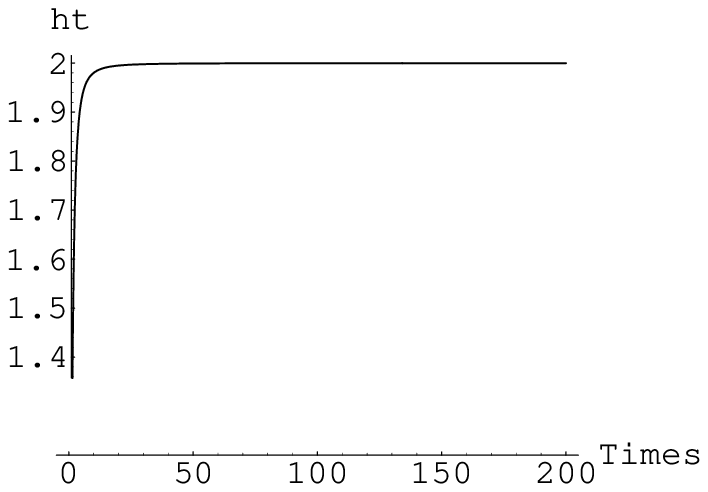}
}
\caption{Representation under alternative A1 of $h_t=2t^2/(1+t^2)$ for time length = 20 (a) and time length = 200 (b)}\label{figA1}
	\end{center}
	\end{figure}

\begin{figure}[H]
	\begin{center}
\subfigure[Time length = 20]
{
		\includegraphics[scale=0.7]{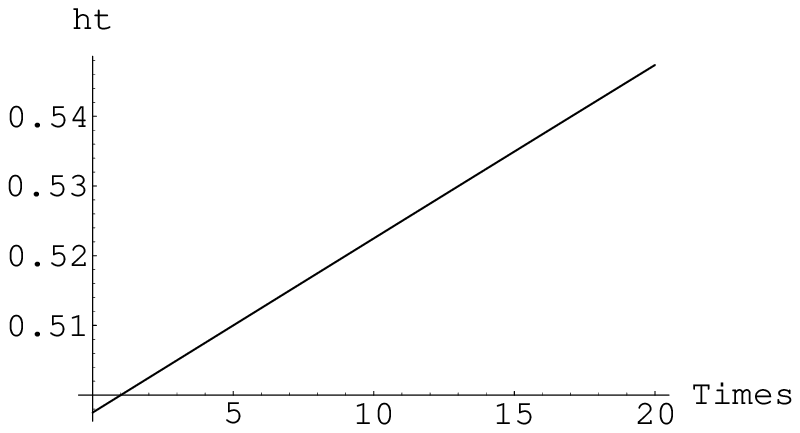}
}
%
\subfigure[Time length = 200]
{		\includegraphics[scale=0.68]{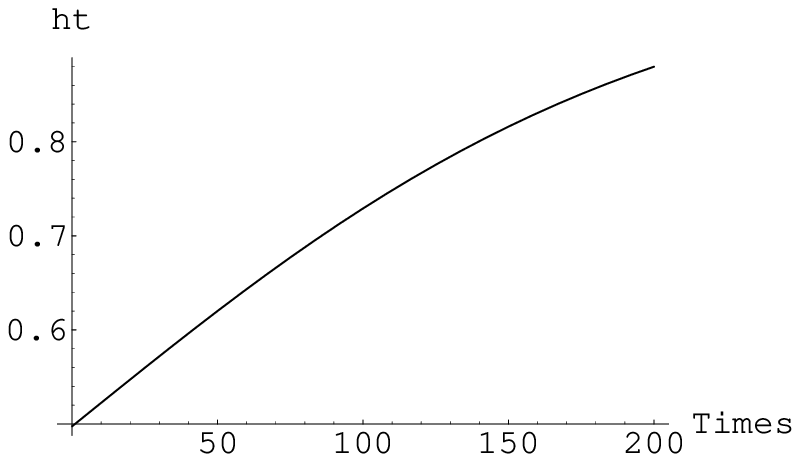}
}
\caption{Representation under alternative A4 of $h_t=(1+\exp(-0.01(t-1)^2))^{-1}$ for time length = 20 (a) and time length = 200 (b)}\label{figA4}
	\end{center}
	\end{figure}

~\\
All alternatives are smooth and are less rough than classical rupture  on the mean or on the variance, except A3 which coincides with a jump on the mean.
The first two alternatives A1-A2 tend quickly to the null model under $H_0$ when the length $n$ increases.
Figure \ref{figA1}  illustrates the proximity of $h_t$ to a constant for large times length in the case of alternative A1. In opposition, alternatives A4-A5 are very smooth and  converge slowly to the null model.  Figure \ref{figA4} illustrates this smooth convergence under alternative A4.

\paragraph{Bootstrap procedure.} To evaluate the power of our
 testing procedure we first consider a Monte Carlo  statistic. Given $M$ points $x_1, \cdots, x_M$ in ${\cal Y}$ we consider
\begin{eqnarray}
S_{M}(w)=\displaystyle\frac{1}{M}\displaystyle\sum_{i=1}^{M}w(x_{i})
A(x_i),
\label{stat}
\end{eqnarray}
where
\begin{eqnarray*}
A(x_i)&=&\max_{1\leq k\leq n}\left\lvert\displaystyle%
\frac{1}{\widehat{\sigma }_{n}(x_{i})\sqrt{n}}\left( \displaystyle%
\sum_{t=1}^{k}\widehat{h}_{t}(x_{i})-{k}\widehat{h}(x_{i})\right)\right\rvert ,
\end{eqnarray*}%
with
\begin{eqnarray*}
\left\{
\begin{array}{lll}
\widehat{\sigma }_{n}^{2}(x) & = & \displaystyle\frac{1}{n}\displaystyle%
\sum_{t=1}^{n}(\widehat{h}_{t}(x)-\bar{h}(x))^{2} \\
\bar{h}(x) & = & \displaystyle\frac{1}{n}\displaystyle\sum_{t=1}^{n}\widehat{%
h}_{t}(x).%
\end{array}%
\right.
\end{eqnarray*}
The convergence 
of the statistic $S_M(w)$  is not guaranteed since the $A(x_i)$ are dependent.
To carry out this problem,  a bootstrap  procedure is proposed.
We construct a naive bootstrap statistic; that is, the test statistic
$S_M(w)$ given in (\ref{stat})
is  compared to the empirical bootstrapped distribution obtained from
$({S_M}^{*b})_{b=1,\cdots,B}$,
 with $S_M^{*b}$  constructed from the bootstraped sample drawn  randomly with replacement and satisfying
 the size equalities
 $N_x^*=N_x$ and $N_y^*=N_y$.
We fix $w$ as a constant.
Note that if $X$ and $Y$ are paired, the bootstrap  procedure consists in drawing randomly with replacement $N$ pairs $(X,Y)$ from the data.  We fix  $B=200$ bootstrap replications.


\paragraph{Powers.}

For each alternative, the test
statistic is computed, based on  sample sizes $N=50, 100, $  for a theoretical level $\alpha=5\%$. The lengths of time's intervals are $n=20, 100 $ and $200$; that is, the function $h_t$ is observed $N$ times for each $t$ varying in $[0;20]$, or  $[0;100]$, or   $[0;200]$, with a step equal to one.
The empirical power of the test is  defined as the percentage of rejection of the
null hypothesis over $10000$ replications of the test statistic under the alternative.

\begin{figure}[H]
	\begin{center}
		\includegraphics[scale=0.5, angle=270]{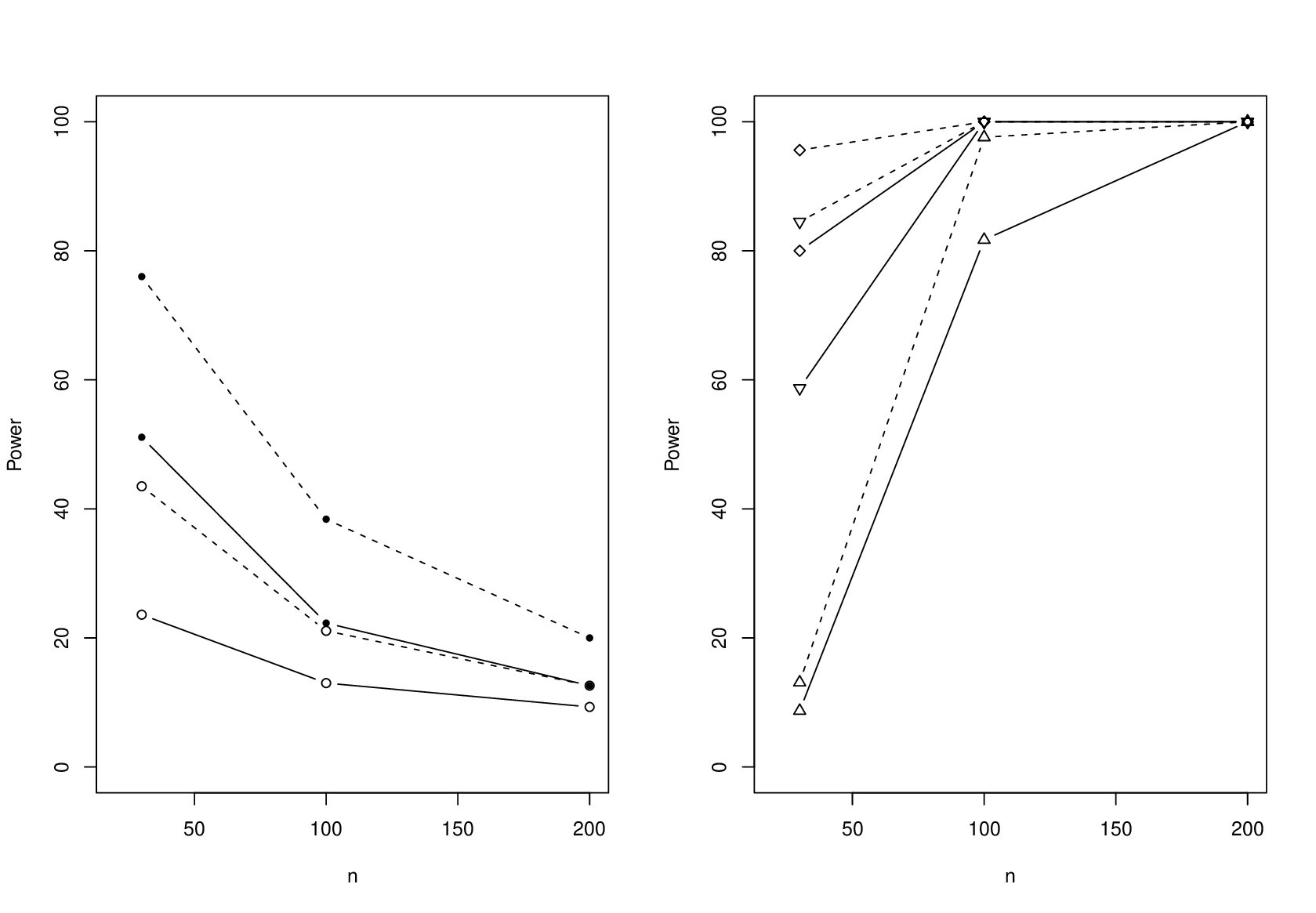}
		\caption{Empirical powers for alternatives A1 ($\bullet$) and A2 ($\circ$) on the left, A3 ($\diamond$), A4 ($\vartriangle$) and A5 ($\triangledown$) on the right, with $X_t$ distributed as ${\cal N}(0,1)$. Solid lines  correspond to $N=50$ and dotted lines correspond to $N=100$. The  lengths of time's intervals are $n=20, 100, 200$}\label{fig1}
	\end{center}
	\end{figure}

Figure \ref{fig1} presents empirical powers of the bootstrap test for all alternatives, in the case where $X_t$ are independent  standard Gaussian variables.  Solid lines and dotted lines correspond to $N=50$ and $100$ respectively.
It can be observed that the power decreases with the length for alternatives A1 and A2. It is in accordance with the previous remark: $h_t$ is close to the null hypothesis for relatively large values of $n$.
Then passing from a time length equal  20 to a time length equal to 200 corresponds to adding variables  with nearly constant transformation in distribution (see Figure \ref{figA1}).

Alternatives A4-A5 have similar behaviors, with a power increasing with $n$. It can be explained by the very slow convergence to the null model.  Here,  passing from a time length equal  20 to a time length equal to 200 corresponds to adding new observations  with a time depending  transformation (see Figure \ref{figA4}).

It is also  observed that power associated to  alternative A3 increases with $n$.

\begin{figure}[H]
	\begin{center}
		\includegraphics[scale=0.5, angle=270]{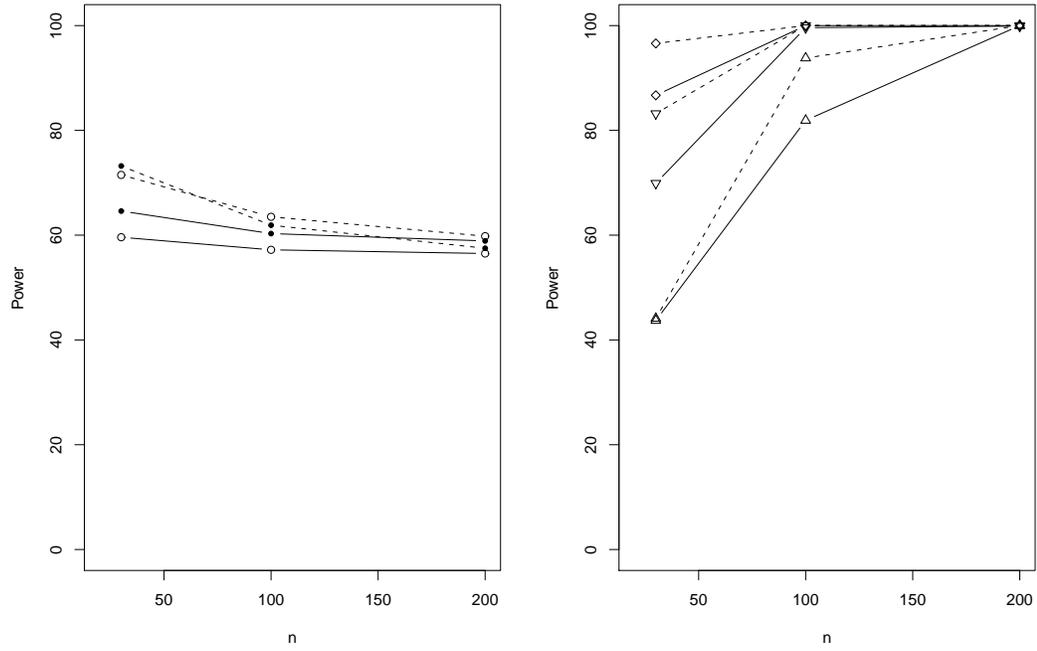}
		\caption{Empirical powers for alternatives A1 ($\bullet$) and A2 ($\circ$) on the left,  A3 ($\diamond$), A4 ($\vartriangle$) and A5 ($\triangledown$) on the right, with $X_t$  following an AR1 process with correlation 0.1. Solid lines  correspond to $N=50$ and dotted lines correspond to $N=100$. The  lengths of time's intervals  are $n=20, 100, 200$}\label{fig2}
	\end{center}
	\end{figure}

In Figure \ref{fig2}  empirical powers are presented in the case where $Y_t$ follows an AR1 process with a correlation coefficient equal to 0.5.
Here powers are slightly better and more stable with respect to the length. This is due to the correlation  inducing  more stability of the process $Y_t$ and permitting a better estimation of $h_t$.

\begin{figure}[H]
	\begin{center}
		\includegraphics[scale=0.5, angle=270]{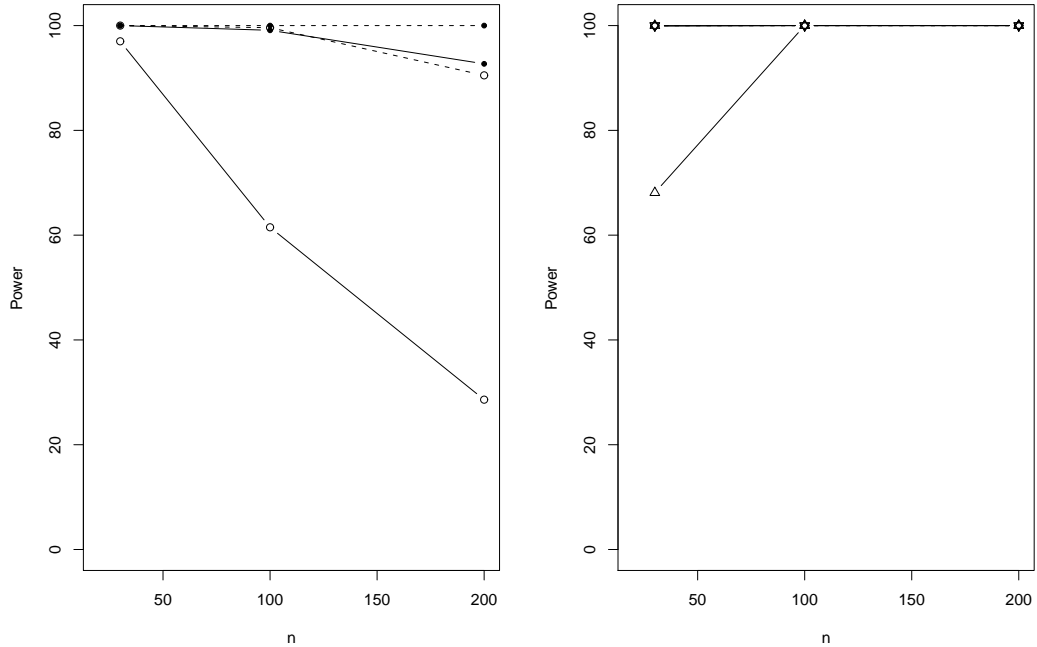}
		\caption{Empirical powers for alternatives A1 ($\bullet$) and A2 ($\circ$) on the left, A3 ($\diamond$), A4 ($\vartriangle$) and A5 ($\triangledown$) on the right, with $X_t$ and $Y_t$ paired. Solid lines  correspond to $N=50$ and dotted lines correspond to $N=100$. The  lengths of time's intervals are $n=20, 100, 200$}\label{fig3}
	\end{center}
	\end{figure}

Figure \ref{fig3} presents results in the case of paired data, with $Y_t$ normally distributed.
Powers are good, due to the fact that transformations occur not randomly since we have considered $X_t=h_t(Y_t)$. Then  $h_t$ can be efficiency estimated and its variations are well detected.





\section{Concluding remarks}

The proposed method concerns the comparison of two processes when panel data are available.
The test permits to detect a change in the relation between the two process distributions. Therefore it can detect
a change in a higher moments (not only in the mean and/or in the variance as almost tests do in this framework).
The asymptotic distribution of the proposed statistic  was derived  under the null of no change in the relation between the two process distributions.

The Monte Carlo simulations show that our test performs well in finite sample and  has a good power against either abrupt or smooth changes.
It is also valid for paired processes and then it can be used to detect a change in $h_t$ in the relation
$X_t=h_t(Y_t)$ (see the paired case in our simulations).
The test  can also be used  as a first step permitting to legitimate
estimation and interpretation of a constant transformation $h$ between two panel data, as for instance
in a medical follow-up study.

 A direction for future research	is to consider a $d$-sample  comparison of distributions, for $d>2$, in the way of \cite{bal8,bal9}. Another direction should consider multivariate distributions.





%
%


\end{document}